\documentclass[11pt]{article}               % The use of LaTeX2e is preferred.
\usepackage{float} %for figures
\usepackage{amssymb,amsfonts,amsmath}
\usepackage{graphicx}              % Include this line if your  document contains figures,
\newcommand {\mat}      [1] {\left[\begin{array}{#1}}
\newcommand {\rix}          {\end{array}\right]}

\title{Numerical methods to compute a minimal realization of a port-Hamiltonian system}                                        % Title,

\author{Karim Cherifi\footnotemark[4] \, and Volker Mehrmann\footnotemark[1] \, and Kamel Hariche\footnotemark[4]
}
\begin{document}

\maketitle

\begin{abstract}                          % Abstract of not more than 200 words.
Port-Hamiltonian (pH) systems are a very important modeling tool in almost all areas of systems and control, in particular in network based model of multi-physics multi-scale systems. They lead to remarkably robust models that can be easily interconnected.  This paper discusses the derivation of pH models  from time-domain input-output data. While a direct construction  of pH models is still an open problem,  we present three different indirect numerical methods for the realization of pH systems. The algorithms are implemented in MATLAB and their performance is illustrated via several numerical examples.
\end{abstract}
\renewcommand{\thefootnote}{\fnsymbol{footnote}}

\footnotetext[4]{
Institute of Electrical and Electronic Engineering, M'hamed Bougara University, Algeria.
\texttt{cherifikarimdz@gmail.com, khariche@yahoo.com}}

\footnotetext[1]{
Institut f\"ur Mathematik MA 4-5,
TU Berlin, Str. des 17. Juni 136,
D-10623 Berlin, FRG.
\texttt{mehrmann@math.tu-berlin.de}. Supported by {\it the German Federal Ministry of Education and Research BMBF within the project EiFer} and by  the German Research Foundation (DFG) via the project {\it Computational Strategies for Distributed Stability Control
in Next-Generation Hybrid Energy Systems}
within Priority Program 1984, Hybrid and Multimodal Energy Systems.}

\section{Introduction}

In this paper we consider the construction of minimal realizations of \emph{linear time-invariant (LTI) port-Hamiltonian (pH) systems} from input-output data in the time domain. Such pH systems are structured representations of dynamical systems \cite{AltS17,BadMBM18,GolSBM03,JacZ12,OrtSMM01,Sch04,Sch06,SchJ14} that typically arise from \emph{energy based modeling} via bond graphs \cite{Bre08,CouJMTB08}. For several application domains this modeling framework has  been implemented in the package, {\sc 20-sim} (\texttt{http://www.20sim.com/}).

The pH framework is extremely powerful, since it encodes underlying physical principles such as conservation laws directly into the structure of the system model. Standard  LTI pH systems take the form
 \begin{eqnarray}
\dot x&=&\left(J-R\right)Qx +(F-P)u,\nonumber \\
y&=&(F+P)^T Qx  + (S+N) u. \label{PHdef}
\end{eqnarray}
Here $J=-J^T \in \mathbb R^{n,n}$, $R=R^T\in \mathbb R^{n,n}$ is positive semidefinite; $F\pm P\in\mathbb{R}^{n,m}$ are \emph{port} matrices; $S+N$, with $S=S^T\in \mathbb R^{m,m}$ and $N=-N^T \in \mathbb R^{m,m}$ is the feed-through from input to output, and $Q=Q^T\in \mathbb R^{n,n}$ positive semidefinite is  associated with a quadratic energy functional \emph{(Hamiltonian)} ${\mathcal H}(x)=\frac 12 x^TQx$. There is an intricate and subtle relationship between pH systems, positive real systems and passive systems, see e.g. \cite{BeaMV19,ByrIW91}. To assure that an LTI system of the form \eqref{PHdef} is Lyapunov stable and passive one further requires that the block matrix
\begin{equation} \label{Kdef}
W=\left[\begin{array}{lc}
R & P \\[1mm]
P^T & S
\end{array}\right]
\end{equation}
associated with the \emph{passivity of the system}, is also positive semidefinite (which we denote as $W\geq 0$ or $W>0$ in the positive definite case). The pH structure has many further spectral properties  \cite{MehMW18}, it is very robust under structured perturbations \cite{MehMS16,MehMS17}, the family of pH systems is closed under power-conserving interconnection \cite{Kle13}, and model reduction of pH systems via Galerkin projection yields (smaller) pH systems \cite{BeaG11,GugPBS12,PolS10}. It also has recently been shown how to robustly transform a passive and Lyapunov stable standard LTI system into a robust pH system \cite{BeaMV19}.

%For many physical domains  first principles based models in pH are available, but
Another important topic that has recently been addressed is the incorporation of constraints in the system models. Then one obtains pH differential-algebraic equations (pHDAEs), or pH descriptor systems, see \cite{BeaMXZ18,Sch13}, which in the LTI case take the form
 \begin{eqnarray}
E \dot x&=&\left(J-R\right)Qx +(F-P)u,\nonumber \\
y&=& (F+P)^T Qx  + (S+N) u. \label{DAEPHdef}
\end{eqnarray}
In the descriptor case it is allowed that $E,Q\in \mathbb R^{n,\ell}$ are rectangular, provided that $Q^TE=E^TQ\geq0$ and it is required that $Q^TJQ=-Q^TJ^TQ \in \mathbb R^{n,n}$, $Q^TRQ=Q^TR^TQ\in \mathbb R^{n,n}$, that the quadratic Hamiltonian is ${\mathcal H}(x)=\frac 12 x^TE^TQx$ and the passivity matrix satisfies
\[
W=\left[\begin{array}{lc}
Q^TRQ & Q^TP \\[1mm]
P^TQ & S
\end{array}\right]\geq 0.
\]
It has been shown that the discussed  properties of standard pH systems carry over to pHDAEs in \cite{BeaMXZ18,HauMM19_ppt,Sch13,SchM18}.

Considering the remarkable system properties and the fact that in many applications a detailed first principle physical modeling of a complicated technical system may be too demanding, it is natural to study  the question of how to construct realizations of pH (or pHDAE) systems from input-output data in the time domain. In particular, it is of great importance  to develop numerical methods for the construction of such realizations of minimal order, usually called minimal realizations. Most current pH modeling approaches work directly on the physical network realization, see e.g.  \cite{SchJ14,XuZZQ12}, while direct data based approaches so far are rare, and direct pH modeling from time-domain data is an open problem.

An example, where a data based approach is urgently required, is the modeling of gas transport networks as pHDAE systems. A general hierarchy of  first principles based models has been derived in~\cite{DomHLT17}. For the flow in the pipe network, after discontinuous Galerkin discretization, see \cite{EggKLMM18}, the model takes the pHDAE form \eqref{DAEPHdef} with $Q=I$, $P=0$, $S+N=0$,
\begin{align} \label{eq:matrices}
E=\begin{bmatrix} M_1 & 0 & 0 \\ 0 & M_2 & 0 \\ 0 & 0 & 0 \end{bmatrix},
\
A=\begin{bmatrix} 0 & G & 0 \\ -G^\top & D & -N^\top \\ 0 & N & 0\end{bmatrix},
\
B=\begin{bmatrix} 0 \\ B_2 \\ 0 \end{bmatrix},\
\end{align}
where $M_1=M_1^T,M_2=M_2^T>0$ are  mass matrices and $x_1$, $x_2$ are coefficient vectors associated with pressure and mass flux, respectively, and $x_3$ represents Lagrange multipliers associated with conservation laws for mass and momentum at network nodes. The matrix $B_2$ is associated with the transport in and out of the pipe network and can be used for coupling with other network components, in particular the compressor
station for which no first principles model, but  huge real world data sets of time domain input-output data are available. In the context of rewriting the whole gas transport system as a network of hierarchies of pH submodels, it would be ideal to generate a model component in pH form.

Another example, where PH modeling from data would be important, is in the modeling of a cable-driven parallel robot (CDPR) \cite{SchYSB18}. PH modeling in this case is used to design physics-shaping controllers, and the PH model was derived analytically. The analysis in \cite{SchYSB18} shows the limitations and the complexity of trying to derive large and complex nonlinear PH models analytically.
%This motivates the development of a data driven PH modeling method.
%The PH model derived has the form
%
%\begin{equation} \label{robot}
% \begin{array}{rcl} \left[ \begin{matrix}
%   {\dot{q}}  \\
%   {\dot{p}}  \\
%\end{matrix} \right] &=& (J-R) \left[ \begin{matrix}
%   \frac{\partial H(p,q)}{\partial q}  \\
%   \frac{\partial H(p,q)}{\partial p}  \\
%\end{matrix} \right] +\bar{G}(q,p)u \\
%  y &=& {{{\bar{G}}}^{T}}(q,p) \left[ \begin{matrix}
%   \frac{\partial H(p,q)}{\partial q}  \\
%   \frac{\partial H(p,q)}{\partial p}  \\
%\end{matrix}  \right]
%\end{array}
%\end{equation}
Similar requirements arise e.g. in complex nonlinear pH modeling for power electronics \cite{CupGBZJMD19} and continuous stirred tank reactors \cite{HoaCJL11}.

In this paper we present three different numerical methods for generating such structured (minimal) pH system realizations. All these require first an interpolation step where the input-output sequence is interpolated to obtain a descriptor (DAE) system.  The resulting realizations from the interpolation step are also usually not well formulated for systems and control applications. For this reason  they are followed by a regularization procedure, see  \cite{BinMMS15_ppt,CamKM12}, which also guarantees that the structural properties are satisfied.
%To derive such an interpolation step that directly generates a system in pHDAE structure is an %\textbf{open problem}.

In the first method, a straightforward approach based on a combination of well-known classical computational control techniques is presented. A model is generated from classical unstructured realization  techniques; this model is converted to a standard minimal state space system and then a linear matrix inequality (LMI) is solved. Using the solution of the LMI, a state space transformation is determined that takes the system to standard pH form.

The second  method also uses a combination of standard techniques but replaces the transformation to a minimal system and the solution of the LMI by a  positive real balanced truncation.

The third method follows a suggestion by \cite{GilS18} that was introduced for dissipative Hamiltonian systems and tackles the minimal pH realization problem as an optimization problem by computing the nearest minimal pH system to a given unstructured realization.

All three approaches are implemented in MATLAB. We compare the three approaches, and illustrate their performance at hand of several numerical examples. %Although these methods deliver reasonably good %results, in principle, we would prefer to have a method that models the pH system directly from %input-output data. However, to obtain such a method is currently an open problem.

\section{Realization and regularization of descriptor systems}\label{sec:interpolation}
In this section we recall the basis for realizing and regularizing a linear time-invariant
descriptor (DAE) system
\begin{equation} \label{DAE}
 \begin{array}{rcl}  E \dot{\hat x} & = & A\hat x + B \hat u, \quad \hat x(0)=0,\\
\hat y&=& C \hat x+D \hat u.
\end{array}
\end{equation}
In a realization process, if the given input-output data are sampled frequency response data, then the Loewner approach \cite{AntLI17,MayA07} is currently one of the most promising and best analyzed method.
If the given data is in terms of input-output time-domain data, the case that we study here, then a modified Loewner approach has recently been proposed in \cite{PehGW17}. This method generates an LTI descriptor system of the form (\ref{DAE}) and is particularly useful if the system is given as a black box routine that simulates the  model from chosen inputs or if the input-output data are just given as a time series.
%The resulting system is an LTI descriptor system of the form (\ref{DAE}) which is constructed directly %from time domain input-output snapshots.
% without prior knowledge of frequency response data. Of course, the more knowledge about the internal %characteristic properties of the system, the more accurate the model will be.

To construct a realization, we assume that inputs $u$ and the resulting outputs $y$ of the system are evaluated or measured at a sequence of time points $t_i$ and recorded as $k$ pairs $(u_i,y_i)$ of vectors
%$k_{\mathrm min} \leq k$,  vectors
$u_i,y_i\in \mathbb R^m$, $i=1,\ldots,k$, where $k$
%_{\mathrm min}$
is selected to ensure that the outputs have typically entered a steady state. Note that for a pH realization we would need the same number of inputs and outputs.
% after $k_{\mathrm min}$ time steps.
The choice of the input sequence $u_i$, the time interpolation points $t_i$,  and the number of samples $k$ are very flexible and can be adapted to best correspond to the real system at hand.
%, including which range of frequencies is of importance.
%Specifically, the range of the possible interpolation frequencies is chosen as
%%
%\begin{equation} \label{rangeF}
%\left[ \frac{2\pi }{k},\frac{2\pi (k-1)}{k} \right]\subset \mathbb{R},
%\end{equation}
%
%so that increasing the number of samples $k$, increases the range of possible frequencies.

The resulting model realization is typically a linear descriptor system of the form (\ref{DAE}),  see \cite{AntLI17,MayA07,PehGW17}.
%used to construct the pH representation. We have used a modified version of the algorithm in %\cite{AntLI17} and transformed the system to a real representation if the system is complex.
%\subsection{Regularization of DAE }\label{sec:regularization}
However, these models are often not so easy to use with current numerical methods for analysis and control design. Often the resulting systems are even close to singular systems, i.e. systems that are not uniquely solvable for a given input $u$. A descriptor system of the form (\ref{DAE}) is called non-regular (singular) if $E,A$ are not square matrices  or in the square case if $\det(sE-A) \equiv 0$.
%If the system is robustly regular (it stays regular even under small perturbations) then no %regularization is necessary.
To check regularity and to obtain a robustly regular system, a regularization procedure as those introduced in \cite{BerV15,BunBMN99,CamKM12,IlcM05a,KunM06}, and implemented for LTI descriptor systems in \cite{BinMMS15_ppt}, needs to  be applied. The regularization procedure of \cite{BinMMS15_ppt,CamKM12} combines a behavioral approach and appropriate feedbacks.
%It possibly uses  renaming of state variables as inputs or inputs as states, removes the feed-through %term by increasing the state-space dimension  and adjusts the dimensions appropriately. %However, no %changes of basis performed.
We will not present this procedure here, see  \cite{BinMMS15_ppt} for a detailed description and a MATLAB implementation. As final output one obtains a reformulated system of the form
\begin{eqnarray}
  {{E}_{1}}\dot{x} &=&{{A}_{1}}x+{{B}_{1}}u, \ x(0)=0, \nonumber \\
   0&=&{{A}_{2}}x+{{B}_{2}}u, \label{eqpar} \\
   y&=&\ Cx +D u, \nonumber
\end{eqnarray}
where
$E_1,A_1 \in \mathbb R^{d,n}$,  $E_2,A_2 \in \mathbb R^{a,n}$, with $n=a+d$, $B_1\in \mathbb R^{d,m}$,
$B_2\in \mathbb R^{a,m}$, $C\in \mathbb R^{m,n}$, $D\in \mathbb R^{m,m}$ and the matrix
$\mat{c} E_1\\ A_2\rix$ is invertible. This property of (\ref{eqpar}) is usually referred to as
the reformulated system being \emph{regular index at most one} \cite{BreCP96} or \emph{strangeness-free} \cite{KunM06} as a free system for $u=0$.

In the form \eqref{eqpar} consistency of initial conditions can be easily checked because these have to satisfy
\begin{equation}
A_{2}x(0)+ B_{2}u(0)=0.
\end{equation}
Note, however, that usually the resulting regularized DAE system will not have a pH structure, so in the next section we will present three procedures to obtain  pHDAE realizations.

\section{PH realization via an unstructured approach} \label{sec:SimpleApp}
In this section we present a straightforward approach to create a pH realization from a regularized  unstructured realization of the form~(\ref{eqpar}). It transforms~(\ref{eqpar}) to a standard state-space form, then to a minimal realization and then this to a pH formulation.

Given a system of the form (\ref{eqpar}), the first step computes matrices $\tilde A,\tilde B, \tilde C, \tilde D$ representing the same input-output map as (\ref{eqpar}) but in standard state space form
\begin{eqnarray}
   \dot{z}&=&\tilde A z+\tilde B u,\ z(0)=0, \nonumber\\
   y&=&\tilde C z+\tilde D u.  \label{StateSp}
\end{eqnarray}
Such a procedure has recently been suggested in \cite{GalKH18}. However, since we have already regularized the descriptor system to the form (\ref{eqpar}), this procedure can be simplified and made more numerically robust. As a first step (since by construction $E_1$ has full row rank, we perform a singular value decomposition, see e.g. \cite{GolV96}, $E_1=U_1 \mat{cc} \Sigma_E & 0 \rix V^T$, with orthogonal matrices $U_1,V$, and $\Sigma_E\in \mathbb R^{d,d}$. We form
\begin{eqnarray*}
\mat{cc} A_{11} & A_{12} \rix&=& U_1^T A_1 V, \ \mat{cc} A_{21} & A_{22} \rix=  A_2 V, \\
\tilde B_1&=& U_1^T B_1,\ \mat{cc} C_1 & C_2 \rix =C V,\ \mat{c} x_1 \\ x_2 \rix= V^T x,
\end{eqnarray*}
partitioned accordingly. As a direct consequence of the regularization procedure
%and the fact that the system (\ref{eqpar}) is regular and of index at most one as a free system with %$u=0$,
it follows that $A_{22}$ is invertible. Multiplying the first block row by $\Sigma_E^{-1}$ and performing a Schur complement with respect to the matrix $A_{22}$, we obtain a standard state-space system of the form \eqref{StateSp} with

\begin{eqnarray*}
\tilde A&=& \Sigma_E^{-1}(A_{11} -A_{12}A_{22}^{-1} A_{21})\in \mathbb R^{d,d},\ \tilde B= \Sigma_E^{-1}( \tilde B_1-A_{12}A_{22}^{-1} B_2)\in \mathbb R^{d,m},\\
\tilde C&=& C_1-C_2A_{22}^{-1} A_{21}\in \mathbb R^{m,d},\ \tilde D=D - C_2 A_{22}^{-1} B_2\in \mathbb R^{m,m},
\end{eqnarray*}
together with a constraint equation $x_2=-A_{22} ^{-1} (A_{21} x_1 + B_2 u)$ which
%restricts the initial conditions and
fixes $x_2$ in terms of $x_1$ and $u$. For further construction we only consider the system in $x_1$, since it represents the transfer function in frequency domain given by
\[
\tilde G(s)= \tilde C (sI-\tilde A)^{-1} \tilde B+ \tilde D.
\]
%
%since the constraint equation has lead to a modification of the feed-through term.
%from $u$ to $y$,
But we should keep in mind that the constraint equation cannot be ignored completely, since it restricts the initial conditions for $x,u$ in the behavior setting. Note that even if the original system did not have a feed-through term, the algebraic constraints will lead to such a feed-through term.
%, see also \cite{CamKM12}.

The presented reduction procedure is computationally acceptable if the condition numbers, see e.g. \cite{GolV96}, of $\Sigma_E$ and $A_{22}$ are small, in relation to the modeling, interpolation, and round-off errors in the regularized system \eqref{eqpar}. If these condition numbers are too large, then this system is close to a singular or high index system and should be treated with great care.

The resulting system in standard state space form \eqref{StateSp} may not be of minimal order. Using the well-known result, see e.g. \cite{Che98}, that a standard state space system of the form (\ref{StateSp}) with a proper frequency domain transfer function $\tilde G(s)$ is minimal if and only if $(\tilde A, \tilde B)$ is controllable and $(\tilde A, \tilde C)$ is observable, one can make the system minimal by removing the uncontrollable and observable parts, e.g. by the techniques
in \cite{Che98,CheH17,DeS00,Kai80,Ros70}, in particular, using the numerically stable staircase form of \cite{Van79} and its implementation in the SLICOT library \cite{BenMSVV99,MehV15}.

Then, keeping only the part that is controllable and observable, results in a minimal realization
\begin{eqnarray}
  \dot {\hat z}\,&=&{\bar A}\hat z+{\bar B} u, \nonumber  \  \hat z(0)=0,
 \\ y&=&{\bar C} \hat z+{\bar D} u, \label{minrealsys}
\end{eqnarray}
with a minimal proper transfer function $\hat G(s) = \hat D+\hat C (sI-\hat A)^{-1} \hat B$. Note that by this minimalization procedure the state-space dimension may have been further reduced.

Each step of this four-stage process of \emph{realization, regularization, transformation to standard state-space form, and making the system minimal} is well understood from its numerical properties, and also concerning its difficulties and robustness under perturbations, thus we have refrained here from presenting more details than necessary.

%For a system that is minimal and also stable and passive it is well-known that it can be transformed by %a state space transformation to a pH system, provided that the number of inputs is equal to the number %of outputs. If this is not the case, one can patch either $\hat B$ with some zero columns and add some  %virtual inputs or add vanishing outputs and output equations to the system to make the numbers equal.
Using the well-known relation between pH systems and systems that are (asymptotically) stable and (strictly) passive, \cite{Wil71,Wil72a}, see also \cite{BeaMV19} for a discussion of some limiting cases, the transformation  to pH form can then be determined by computing a positive definite solution $\hat Q$ to the linear matrix inequality (LMI)
	\begin{equation} \label{LMI}
	\left[ \begin{matrix}
   {{\hat A}^{T}}\hat Q+\hat Q \hat A & \hat Q\hat B-{{\hat C}^{T}}  \\
   {{\hat B}^{T}}\hat Q-\hat C & -(\hat D+{{\hat D}^{T}})  \\
\end{matrix} \right]\le 0.
\end{equation}
See also \cite{BeaMX15_ppt} for a more detailed discussion on  numerical properties and a direct computation of this solution via a structured eigenvalue problem.

Note the LMI (\ref{LMI}) in general may not have a positive definite solution if $D+D^T$ is singular. A commonly used practice is to use a perturbation $\hat D +\hat D^T+\varepsilon I$ to make this term positive definite. However, there is no perturbation analysis of the effect of this  perturbation on the solution procedure as a whole. Another, more reliable approach is to first perform a projection of the problem, by computing an orthogonal  matrix
$U_2$ such that
\[
U_2^T (D+D^T) U_2=\mat{cc} \hat D_1+ \hat D_1^T & 0 \\ 0 & 0 \rix,\ \hat BU_2=\mat{cc} \hat B_1 & \hat B_2 \rix,\ U_2^T \hat C=\mat{c} \hat C_1 \\ \hat C_2 \rix,
\]
where $\hat D_1+\hat D_1^T>0$. Then one has to compute a positive definite solution of the LMI
	\begin{equation} \label{LMIred}
	\left[ \begin{matrix}
   {{\hat A}^{T}}\hat Q+\hat Q \hat A & \hat Q\hat B_1-{{\hat C_1}^{T}}  \\
   {{\hat B_1}^{T}}\hat Q-\hat C_1 & -(\hat D_1+{{\hat D_1}^{T}})  \\
\end{matrix} \right]\le 0
\end{equation}
satisfying the constraint  $\hat Q \hat B_2=\hat C_2^T$. According to the results of \cite{Wil71,Wil72a},
if the system is stable and passive then such a $Q=Q^T>0$ exists. If not, then a further regularization procedure is needed that can be based on \cite{BeaMX15_ppt}, however, in the following we assume that such a $\hat Q>0$ can be computed. Then one can determine a Cholesky factorization  $\hat Q=T{T}^T$, with a nonsingular matrix $T$. This Cholesky factor can be used to transform  \eqref{minrealsys} to pH form as in \cite{BeaMV19,BeaMX15_ppt} leading to the system matrices
\begin{eqnarray}
   Q&=&I,\
 J=\frac{1}{2}(TAT^{-1}-(TAT^{-1})^T)\ R=-\frac{1}{2}(TAT^{-1}+(TAT^{-1})^T), \nonumber \\
 F&=&\frac{1}{2}(TB+{{(C{{T}^{-1}})}^{T}}),\ P=\frac{1}{2}(-TB+{{(C{{T}^{-1}})}^{T}}),\label{solPH}\\
 S&=&\frac{1}{2}(D+{{D}^{T}}),\ N=\frac{1}{2}(D-{{D}^{T}}).\nonumber
\end{eqnarray}

\begin{figure}
\centering
  \includegraphics [width=0.7\linewidth]{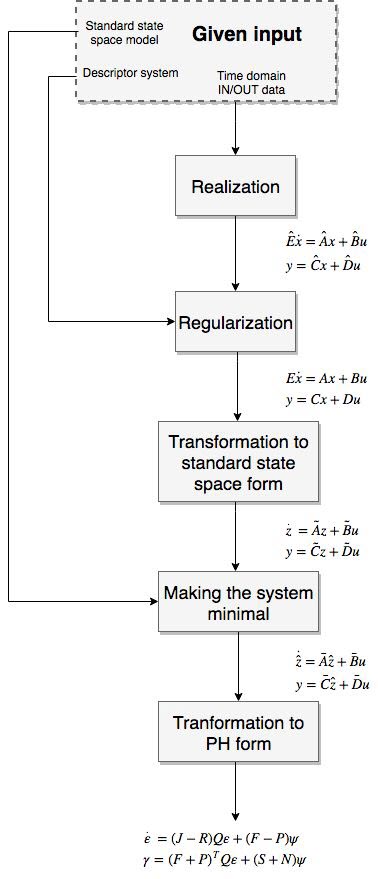}
  \caption{Block diagram of the simple unstructured approach}
  \label{fig:Simple}
\end{figure}

We have implemented this procedure in MATLAB using the LMI approach as well as a Riccati equation approach \cite{BeaMV19}, see Figure \ref{fig:Simple} for a block diagram of the procedure.
To avoid unnecessary computation, we allow three different types of starting data for the resulting MATLAB function, either input/output data, or a general descriptor system or a standard state space system. Depending on these data some blocks of the procedure may be skipped.
Numerical results obtained with this procedure are presented in Section~\ref{sec:numerical}.

In this section we have discussed a five step procedure \emph{realization, regularization, transformation to standard state-space form, making the system minimal, transformation to PH form} to generate a standard state-space pH realization from given time domain input-output data. It should be noted that although every single step in the discussed procedure is well understood from the point of view of error analysis and perturbation theory, the procedure may face a significant number of difficulties. These include a singularity or close to singularity of the system \eqref{DAE}, an ill-conditioning of the matrices $\Sigma_E$, $A_{22}$, a close to non-controllability or non-observability of the minimal system \eqref{minrealsys}, or that this system is not stable or not passive. If the system is not stable,  then this can be fixed  by an appropriate feedback, since the system is controllable. If the system is not passive, then a perturbation to the nearest passive system
can be performed. But even if the system is passive it may happen that the Cholesky factor $T$ is ill-conditioned.

It is almost impossible to carry out a detailed error and perturbation analysis for the complete procedure, since there are too many factors that can go wrong in the procedure.
%The five steps all may be sensitive to perturbations such as data, model and numerical errors, and they %require difficult numerical procedures such as rank decisions in finite precision arithmetic, or the %solution of linear systems which may be ill-conditioned.
Furthermore, since each of the steps  builds upon the next, this simple  straightforward approach will not always yield satisfactory results. Finally, as we will see below, it is computationally demanding and does not scale well to large scale problems.

\section{Transformation based on positive real balanced truncation}\label{sec:PRBT}

Another straightforward approach to generate a minimal pH realization from an unstructured standard realization  is to exploit the relationship between pH systems and positive real systems. Consider a DAE system of the form (\ref{DAE}) with a transfer function $G(s)=D+C(sE-A)^{-1} B=p(s) +r(s)$, where $p=\sum_{i=0}^{\nu} M_i s^i$ with ${{M}_{i}}\in {{\mathbb{R}}^{m\times m}}$, is the polynomial and $r$ the proper rational part of $G$.
%,  and assume that
%the system is already strangeness-free (index one) in the behavior as a free system with $u=0$, which %means that $p(s)=D$.
Let $P_r$ and $P_l$ be the spectral projectors onto
the right and left deflating subspaces of $\lambda E-A$ corresponding to the finite eigenvalues.
% and let ${{M}_{k}}\in {{\mathbb{R}}^{m\times m}}$ denote the coefficients of the polynomial part of the system.
It has been shown in \cite{BerR13} that if the system (\ref{DAE}) is \emph{$R$-minimal}, i.e. the dynamical part of the system is controllable and observable and if the transfer function $G$ is passive, then the \emph{projected Lur'e equation}
\begin{align} \label{Lure}
   AX{{E}^{T}}+EX{{A}^{T}} &=-{{K}_{C}}{{K}_{C}}^{T},~~~~X={{P}_{r}}X{{P}_{r}}^{T}\ge 0, \\
  EX{{C}^{T}}-{{P}_{l}}B &=-{{K}_{C}}{{J}_{C}}^{T},~~~~{{M}_{0}}+{{M}_{0}}^{T}={{J}_{C}}{{J}_{C}}^{T} \nonumber
\end{align}
is solvable for ${{K}_{C}}\in {{\mathbb{R}}^{n,m}},{{J}_{c}}\in {{\mathbb{R}}^{m,m}}$ and $X\in {{\mathbb{R}}^{n,n}}$.
Furthermore, if ${M}_{k}=0$ for $k>1$, ${M}_{1}={M}_{1}^{T}\ge 0$ and the
projected Lur’e equations have a solution, then $G$ is positive real.

When, starting from a system of the form (\ref{DAE}), we have  performed the regularization procedure and obtained the system (\ref{eqpar}), then we have ${{M}_{0}}=D$ and $M_\ell=0$ for $\ell>0$. Then, if the resulting system is passive, the solution of \eqref{Lure} can be computed, and also instead of \eqref{Lure}, its dual form can be solved equivalently.
\begin{align} \label{Lure2}
   {{A}^{T}}YE+{{E}^{T}}YA &=-{{K}_{O}}^{T}{{K}_{O}},~~~~Y={{P}_{l}}^{T}Y{{P}_{l}}\ge 0, \\
  {{E}^{T}}YB-{{P}_{r}}^{T}C^{T} &=-{{K}_{O}}^{T}{{J}_{O}},~~~~{{M}_{0}}+{{M}_{0}}^{T}={{J}_{O}}^{T}{{J}_{O}} \nonumber
\end{align}
  for ${{K}_{O}}\in {{\mathbb{R}}^{m,n}},{{J}_{O}}\in {{\mathbb{R}}^{m,m}}$ and $Y\in {{\mathbb{R}}^{n,n}}$ .

The R-minimality condition in the first
part of Theorem 3.1 can be weakened. For standard
systems, less restrictive conditions on the involved matrices are given in \cite{CleALM97,LanR95}.

The solution of the Lur'e equations can then be used to compute the balanced truncation, where
the minimal solutions of the projected Lur'e equations, ${{X}_{\min }}$ and ${{Y}_{\min }}$, respectively, are the positive real controllability and observability Gramians. The system is called positive real balanced if ${{X}_{\min }}={{Y}_{\min }}=\mbox{\rm diag}(\Pi ,0)$ with $\Pi =\mbox{\rm diag} ({{\pi }_{1}},...,{{\pi }_{{{n}_{f}}}})$, where the ${{\pi }_{j}}$ are the positive real characteristic values of the system, i.e.  the square roots of the non-zero eigenvalues of the matrix ${{X}_{\min }}{{E}^{T}}{{Y}_{\min }}E$.

\begin{figure}
\centering
  \includegraphics[width=0.6\linewidth]{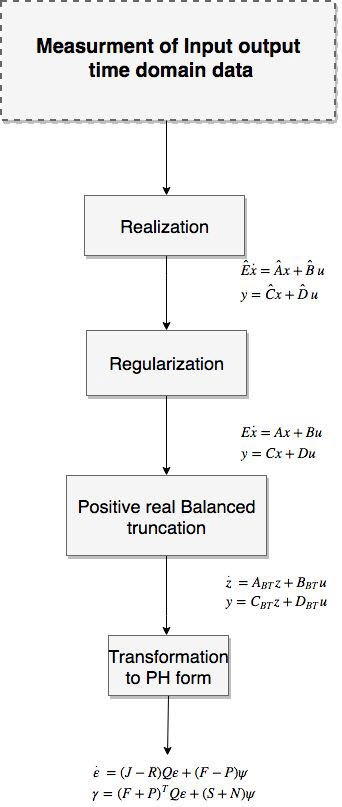}
  \caption{Block diagram of the Positive Real Balanced Truncation Method}
  \label{fig:PRBT}
\end{figure}

After computing the balancing transformation, one can compute  a transformation
\begin{equation} \label{transPRBT}
\left[ \hat{E},\hat{A},\hat{B},\hat{C},\hat{D} \right]=\left[ {{W}_{b}}E{{T}_{b}},{{W}_{b}}A{{T}_{b}},{{W}_{b}}B,C{{T}_{b}},D \right]
\end{equation}
such that a balanced truncation is obtained and states corresponding to insignificant ${{\pi }_{j}}$ are truncated, see \cite{ReiS10}.
If $D+{{D}^{T}}$ is nonsingular, then the positive real Lur'e system is equivalent to the projected Riccati equation
\begin{align}
  & {{A}^{T}}YE+{{E}^{T}}YA+{{({{B}^{T}}YE-C{{P}_{r}})}^{T}}{{(D+{{D}^{T}})}^{-1}}({{B}^{T}}YE-C{{P}_{r}})=0,
  \nonumber \\
 & Y=P_{l}^{T}Y{{P}_{l}}. \nonumber
\end{align}
If $D+{{D}^{T}}$ is rank deficient, then actually the Riccati equation is not well defined and as in the previous section no positive definite solution exists that can be used for the transformation. Again commonly a perturbation ($\varepsilon I$ is used to make $D+{{D}^{T}}+ \varepsilon I>0$. This method has been implemented in the MORLAB toolbox \cite{BenW18}, where the truncation method is applied using Newton-type methods, where Lyapunov equations are solved in each step using the matrix sign function method. For illustration, the Positive Real Balanced truncation (PRBT) method is presented in Figure~\ref{fig:PRBT}.

Again, although there are much fewer steps in this procedure, it is very difficult to carry out a detailed error and perturbation analysis for the complete procedure, since still many steps can go wrong in this procedure. However, the procedure scales much better than that of the previous section for large scale problems.

%The solution $\hat Q$ of the LMI \eqref{LMI} is needed to compute the matrices of the Port Hamiltonian %system  $\hat Q$ in this case can be simply constructed as a diagonal matrix where its entries are the %first r Hankel singular values where r is the order of the resulting system. The different PH system %matrices are then constructed as in \eqref{solPH}.

\section{Nearest Port Hamiltonian system}\label{sec:NearestPH}

Having discussed in the previous section two procedures for the construction of a pHODE realization of a descriptor system arising from a given input-output sequence and having also remarked on the potential numerical difficulties  in computing the  desired pHODE realization, in this section we proceed in a completely different way.  Many of the difficulties discussed in the previous two approaches arise from the fact that the interpolation procedure may lead to a system being a (close to) singular, unstable or non-passive descriptor system even though the system from which the data were generated was originally stable and passive. Thus we may just ask the question whether we can just solve an optimization problem to compute the nearest stable, passive and regular and index at most one (as a free system) phDAE  system. Let us denote the set of such systems by $\mathcal{D}$. For a system with given realization $(E,A,B,C,D)$ we would then have to solve  the minimization problem
\begin{equation}
\underset{(\tilde{E},\tilde{A},\tilde{B},\tilde{C},\tilde{D})\in \mathcal{D}}{\mathop{\inf }}\,\mathcal{F}(\tilde{E},\tilde{A},\tilde{B},\tilde{C},\tilde{D})
\end{equation}
where
\begin{multline*}
\mathcal{F}(\tilde{E},\tilde{A},\tilde{B},\tilde{C},\tilde{D})=\left\| A-\tilde{A} \right\|_{F}^{2}+\left\| B-\tilde{B} \right\|_{F}^{2}+\left\| C-\tilde{C} \right\|_{F}^{2} \\ +\left\| D-\tilde{D} \right\|_{F}^{2}+\left\| E-\tilde{E} \right\|_{F}^{2}.
\end{multline*}
This a very difficult optimization problem, since the feasible set $\mathcal D$ is unbounded,
highly non-convex, and neither open nor closed, see \cite{GilMS18,GilS18} for a detailed analysis.  However, requiring the structure of pHDAE form comes to help, because we can reformulate the problem
as the optimization problem
\begin{equation*}
\underset{M,J,R,F,P,S}{\mathop{\inf }}\,\mathcal{G}(M,J,R,F,P,S)\text{ subject to }{{J}^{T}}=-J,\ M\ge 0,\ \left[ \begin{matrix}
   R & P  \\
   {{P}^{T}} & S  \\
\end{matrix} \right]\ge 0,\
\end{equation*}
where
\begin{eqnarray*}
\mathcal{G}(M,J,R,F,P,S)&=&\left \| E-M \right \|^2 + \left\| A-(J-R) \right\|_{F}^{2}+\left\| B-(F-P) \right\|_{F}^{2} \\ && +\left\| C-{{(F+P)}^{T}} \right\|_{F}^{2}
 + \left\|  \frac{D+{{D}^{T}}}{2}-S \right\|_{F}^{2}.
\end{eqnarray*}
For a slightly different version of this optimization problem, see \cite{GilMS18,GilS18}, where
an  algorithm to solve this optimization problem is based on a fast projected gradient method (FGM) which has the advantage of being in general much faster than the standard projected gradient method even in the nonconvex case. The algorithm presented in \cite{GilS18} also includes a restarting procedure which is necessary in our case, since the objective function is not convex and hence FGM without restart is not guaranteed to converge \cite{GhaL16}.  It is also possible to add weights to each term of the objective function.

Choosing good initial points is crucial to obtain good solutions for nonconvex optimization problems and this is still an open problem. In \cite{GilS18}, an LMI based initialization procedure as discussed in our first approach is proposed. This LMI-based initialization works well when the initial system is close to being passive. %Otherwise it may provide rather bad initial points.
However, solving the LMI may require a large computational effort if the system is large scale.
%It is worth noting that the solution found by this algorithm may not be optimal as it may get stuck in a %local minimum, see also \cite{GilMS18}.

\section{Numerical examples}\label{sec:numerical}
In this section we present some numerical examples to illustrate the performance of the different methods. All tests were performed on a laptop Intel Core i7-5600U CPU @2.6GHz, 16GB RAM. The computation times of the algorithms are measured in seconds.

%\subsection{Illustrative example}\label{sec:Simpleex}
Consider  a descriptor system (\ref{DAE}) of order $5$ with
\begin{align}\nonumber
  & E=\begin{bmatrix}
   0 & 0 & 19 & 15 & 5  \\
   0 & 4 & 14 & 13 & 14  \\
   0 & 9 & 10 & 1 & 11  \\
   0 & 7 & 9 & 6 & 12  \\
   0 & 8 & 1 & 17 & 20  \\
\end{bmatrix}, \text{  }A=\begin{bmatrix}
   17 & 10 & 10 & 15 & 7  \\
   9 & 2 & 4 & 6 & 9  \\
   18 & 8 & 20 & 12 & 15  \\
   5 & 1 & 4 & 2 & 19  \\
   14 & 15 & 3 & 3 & 12  \\
\end{bmatrix}, \text{  }B=\begin{bmatrix}
   2  \\
   20  \\
   1  \\
   2  \\
   18  \\
\end{bmatrix}, \text{ } \\ \nonumber
 & \text{  }C= \begin{bmatrix}
   16 & 19 & 3 & 14 & 14  \\
\end{bmatrix}, \text{    }D=9.3.  \nonumber
\end{align}

Applying the first method results in a minimal pHODE system with matrices $E=Q=I$, $S=N=0$,
\begin{eqnarray*}
   J&=&\begin{bmatrix}
   0 & 0.5695 & 0.3115 & -0.5990  \\
   -0.5695 & 0 & -0.1870 & -0.0735  \\
   -0.3115 & 0.1870 & 0 & 1.3185  \\
   0.5990 & 0.0735 & -1.3185 & 0  \\
\end{bmatrix},\ F=\begin{bmatrix}
   6.1006  \\
   5.1569  \\
   -1.3935  \\
   -2.1381  \\
\end{bmatrix}, \\  \nonumber
  R&=&\begin{bmatrix}
   0.5644 & 0.0629 & -0.3760 & 0.1783  \\
   0.0629 & 0.0070 & -0.0419 & 0.0199  \\
   -0.3760 & -0.0419 & 0.2505 & -0.1188  \\
   0.1783 & 0.0199 & -0.1188 & 0.0563 \end{bmatrix},\
   P=\begin{bmatrix}
   0.3527  \\
   0.0393  \\
   -0.2350  \\
   0.1114  \\
\end{bmatrix}
\end{eqnarray*}
of order $4$.

Applying the PRBT method results in a pHODE system with $E=Q=I$, $N=0$, $S=0.2204$,
\begin{eqnarray*}
 J&=&\begin{bmatrix}
   0 & 0.1770 & 0.0321 & -0.9770  \\
   -0.1770 & 0 & -0.7213 & -0.0081  \\
   -0.0321 & 0.7213 & 0 & -1.0257  \\
   -0.9770 & 0.0081 & 1.0257 & 0  \\
\end{bmatrix},\ F=\begin{bmatrix}
   8.37611  \\
   -0.0154  \\
   -0.1999  \\
   -0.1481  \\
\end{bmatrix}, \\
R&=&\begin{bmatrix}
   0.3662 & -0.0566 & -0.0948 & -0.4187  \\
   -0.0566 & 0.0088 & -0.0147 & 0.0648  \\
   -0.0948 & -0.0147 & 0.0245 & -0.1084  \\
   -0.4187 & 0.0648 & -0.1084 & 0.4787  \\
\end{bmatrix},\  P \begin{bmatrix}
   0.2841  \\
   -0.0439  \\
   -0.0736  \\
   -0.3248 \\
\end{bmatrix}.
\end{eqnarray*}
Finally, using the nearest pHDAE algorithm results in a pHODE system with
$S=0.2204$, $N=0$,
\begin{eqnarray*}
Q&=&\begin{bmatrix}
   4.6335 & 4.1536 & 0.7786 & -0.9624  \\
   4.1589 & 10.9526 & -1.1426 & -2.1357  \\
   0.7811 & -1.1365 & 0.9656 & 0.2849  \\
   -0.9665 & -2.1421 & 0.2866 & 0.5735  \\
\end{bmatrix},\  F=\begin{bmatrix}
   0.4406  \\
   0.4741  \\
   0.4523  \\
   -8.6106  \\
\end{bmatrix} \\
 J&=&\begin{bmatrix}
   0 & -0.4905 & 0.2285 & -3.3405  \\
   0.4905 & 0 & -1.3097 & 1.6185  \\
   -0.2285 & 1.3097 & 0 & 6.6314  \\
   3.3405 & -1.6185 & -6.6314 & 0  \\
\end{bmatrix}, \ P=\begin{bmatrix}
   0.3025  \\
   -0.0583  \\
   -0.3741  \\
   0.1620  \\
\end{bmatrix} \\
 R&=&\begin{bmatrix}
   0.4799 & -0.0715 & -0.6606 & 0.5261  \\
   -0.0715 & 0.0206 & 0.0777 & 0.0164  \\
   -0.6606 & 0.0777 & 0.9715 & -0.9760 \\
   0.5261 & 0.0164 & -0.9760 & 1.6382  \\
\end{bmatrix}.
\end{eqnarray*}

Note that $Q=E^{-1}$ is not the identity for the nearest pHDAE algorithm. This constitutes an alternative to the previous methods, but the resulting system has order equal to the original system, and it is necessary to employ a structure preserving minimalization procedure if one wants a minimal pHDAE system.

\begin{figure} [H]
\centering
  \includegraphics[width=\linewidth]{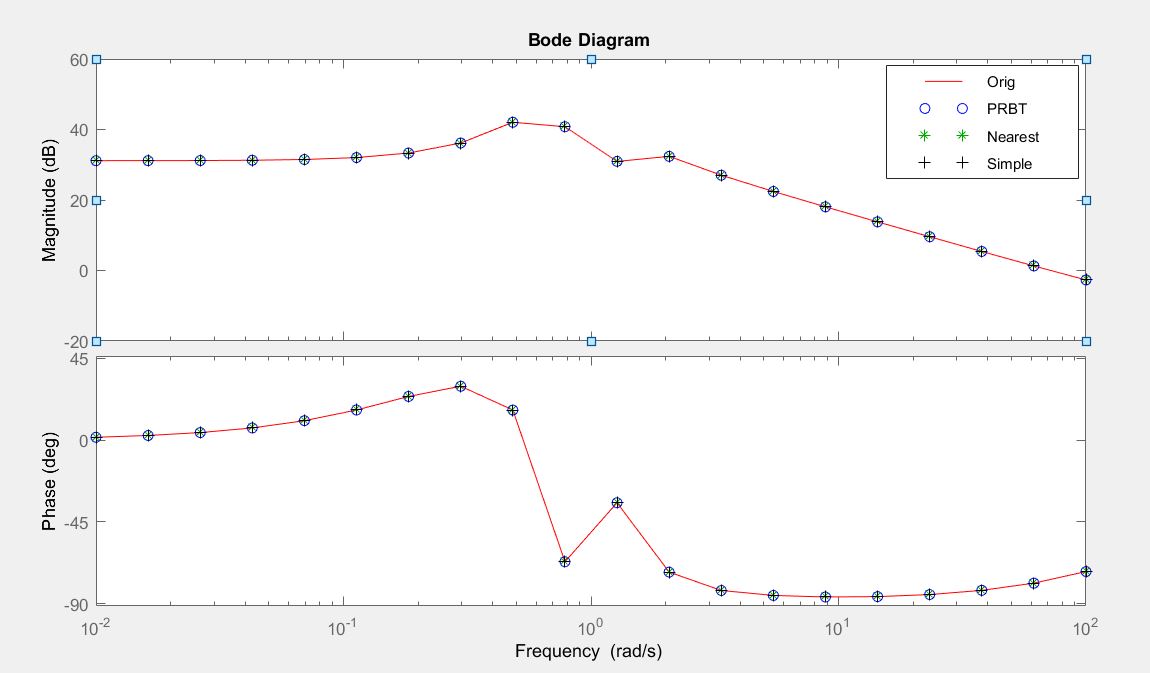}
  \caption{Bode plot comparison for the illustrative example}
  \label{fig:bode1}
\end{figure}

A comparison of the transfer functions of the three pH realizations  with the original system  is shown in the Bode plots in Figure \ref{fig:bode1}.

\begin{figure} [H]
\centering
  \includegraphics[width=\linewidth]{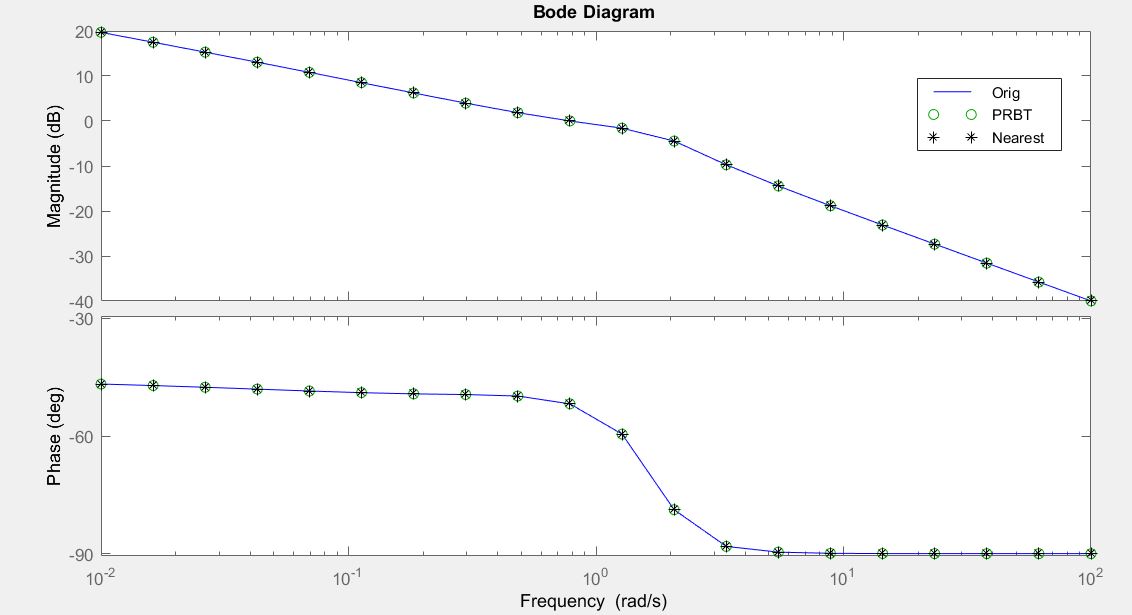}
  \caption{Bode plot comparison for the ladder network example}
  \label{fig:bode2}
\end{figure}

%\subsection{Ladder network }\label{sec:ELECex}
As second example we consider an electrical ladder network consisting of resistors, capacitors and inductors. This example has order $2000$. It is tested for all the presented algorithms. The computation time for the straightforward approach exceeds the allowed time of $3600 s$, this is why it not included in the bode plot comparison in Figure~\ref{fig:bode2}. The resulting pHODE system is obtained without losing the frequency response of the system. Also, in the case of the PRBT approach, the resulting system has order $22$ which is the numerically minimal order. This is a huge improvement compared to the given system which has order $2000$.

%\subsection{Discussion }\label{sec:discussion}
In terms of computation time, the difference between the different methods are large. Table~\ref{tab:table} shows the computation times for the previous examples in seconds $(s)$.

\begin{table}[h!]
  \begin{center}
    \caption{Computation time of the algorithms for the presented examples}
    \label{tab:table}
    \begin{tabular}{l|c|c|c}
    \hline
      \textbf{Examples} & \textbf{Simple method}& 				\textbf{PRBT} & \textbf{Nearest} \\
      \hline
      Illustrative example & 0.82441 & 3.45088 		& 2.56567\\
      \hline
      Ladder network & Failed (time$>3600s$) & 719.74701 & 10.26562\\
%      \hline
%      Mass spring damper & Failed(time$>$3600) & 117.18015 & 14.95489\\
      \hline
    \end{tabular}
  \end{center}
\end{table}

The straightforward simple method presented very good results for the simple illustrative example of order $5$ but does not scale well for larger systems. The nearest pH algorithm gave the best  results in terms of computation time but suffers from the problems mentioned before as it depends on the initialization and it can get stuck in a local minimum. Also it does not lead to a minimal system. The positive real balanced truncation method, although requiring more time than the nearest pH is very robust, and also results in a reduced model.

\section{Conclusion }\label{sec:Conclusion}

In this paper, we presented three methods to obtain a port-Hamiltonian system from time domain input-output data that are realized as a  descriptor system or a standard state space representation. The first approach is a five step procedure that may have many numerical difficulties, is hard to analyze, and does not scale well in computation time. The second method is based on positive real balanced truncation and finally the third  method solves an optimization problem to compute the nearest port-Hamiltonian system. The three procedures were implemented in MATLAB. This implementation and testing indicates that the method based on the positive real balanced truncation is the most robust especially for large scale systems. However, when initiated  well, the nearest port-Hamiltonian system algorithm is the most promising but needs further implementation improvements. It should in particular be combined with a model reduction step as in \cite{HauMM19_ppt} and a minimalization step.

Ideally one would like to construct a pH representation directly from the given input-output data, but this is still an open problem.

%\bibliographystyle{plain}
%\bibliography{CheM}
\end{document}